\documentclass[reqno,b5paper]{amsart}
\usepackage{latexsym,amsmath,amsfonts,amssymb}
\usepackage{graphics}
\usepackage{epsfig}
\usepackage{amsmath,amscd}
\usepackage[all,cmtip]{xy}

\newtheorem{thm}{Theorem}[section]
\newtheorem{deff}[thm]{Definition}

\newtheorem{prop}[thm]{Proposition}
\newtheorem{cor}[thm]{Corollary}

\newcommand{\vG}{\varGamma}

\newcommand{\wt}{\widetilde}

\def\R{{\mathbb R}}

\def\be{\begin{equation}}

\def\ee{\end{equation}}
\def\ba*{\begin{eqnarray*}}

\def\ea*{\end{eqnarray*}}

\newcommand {\f}{\hfill$\Box$}
\newcommand {\Pf}{\bf Proof.}
\newcommand{\bee}{\begin{equation*}}
\newcommand{\eee}{\end{equation*}}

\begin{document}

\title[Pettis integrability of fuzzy]{Pettis integrability of fuzzy mappings with values in arbitrary Banach spaces}
\author[L. DI PIAZZA and V. MARRAFFA ]{L. DI PIAZZA  and V. MARRAFFA}
\thanks{The  authors were partially supported
 by the grant of GNAMPA prot. U 2016/000386 }

\newcommand{\acr}{\newline\indent}
\address{\llap{*\,}Department of Mathematics and Computer Sciences \acr University of
Palermo\acr via Archirafi 34\acr 90123 Palermo\acr ITALY}
\email{luisa.dipiazza@unipa.it}

\address{\llap{*\,}Department of Mathematics and Computer Sciences\acr University of
Palermo\acr via Archirafi 34\acr 90123 Palermo\acr ITALY}
\email{valeria.marraffa@unipa.it}


\bigskip

\begin{abstract}
\noindent
In this paper we study  the Pettis integral of fuzzy mappings in  arbitrary Banach spaces. We present  some properties of the Pettis integral of fuzzy mappings and we give conditions under which a scalarly integrable fuzzy mapping is Pettis integrable.
\end{abstract}

\maketitle

\begin{center}
{\it Cordially dedicated to Professor Paolo de Lucia}\end{center}
\begin{center}
{\it on the occasion of his 80-th birthday with esteem and admiration}
\end{center}

\section{Introduction}

Due to interest in applications (e.g. control theory, optimization or mathematical economics), a wide theory was developed for set-valued integrability in separable Banach spaces (see  \cite{BS}, \cite{castaing valadier}, \cite{hp} and reference inside, \cite{dp1,dp2},
\cite{eh}, \cite{amart}). 
Recently several authors studied set-valued integrability without separability assumptions in
Banach spaces (see  \cite{b1},  \cite{b3}, \cite{b4} \cite{ckr,ckr1}, \cite{cr}, \cite{dp3} and \cite{mu}).

\noindent It is well known the useful role played by  fuzzy mappings in the theory  of fuzzy sets and their applications. Integrals of fuzzy mappings are a  natural extension of  set-valued integrability.  For  fuzzy mappings in ${\mathbb{R}}^n$ an  integral   was introduced  in \cite{Ka} and \cite{Ma}, while the case of  fuzzy mappings in separable Banach space  was considered and studied  in \cite{Xue1}, \cite{Xue} and \cite{Xue2}.  The fuzzy Henstock integral of fuzzy mappings  in  ${\mathbb{R}}$, or  ${\mathbb{R}}^n$, has been introduced and studied in \cite{wg1, wg2} and  \cite{DPfuzzy}, while  in \cite{mufuzzy} the more general case of Banach space valued  fuzzy Henstock integral is considered.

\noindent  Park in \cite{Park1} introduced  the Pettis integral for  fuzzy mappings in separable Banach spaces.

In this paper we continue the study of the Pettis integral for fuzzy mappings in not necessarily separable Banach spaces.  An  important technical tool for our investigation is the following decomposition result (Theorem \ref{tdec})

\vspace{2ex}
\noindent   {\it  $\wt{\vG}\colon \Omega \to {\mathcal F}(X)$ (here  ${\mathcal F}(X)$ is the generalized fuzzy number space associated to the Banach space $X$) is a Pettis integrable fuzzy mapping if and only if
$\wt{\vG}$ can be represented as $\wt{\vG}(\omega)=\wt{G}(\omega)+\wt{f}(\omega)$, where $\wt{G}\colon \Omega \to {\mathcal F}(X)$ is a Pettis integrable fuzzy mapping whose support functions are  nonnegative and
 $\wt{f}$ is  a  Pettis integrable fuzzy mapping generated by a Pettis integrable selection of  $\wt{\vG}$.}

 \vspace{2ex}

\noindent In \cite{mu} (pretty level-wise) results regarding the Pettis integrability of  scalarly integrable multifunction are proved. Here it is shown a kind of  their  ``consistency'' with respect to the fuzzy nature of the fuzzy Pettis integral. Indeed, we  prove that the primitive of a Pettis integrable fuzzy mapping is a generalized fuzzy number measure (Theorem \ref{p1}) and we give a characterization of the Pettis integrable fuzzy mappings by means the notion of core (Theorem \ref{tcore}).

It is well known  that if $X$ is a separable Banach space and $c_0 \not\subseteq X$, then each single $X$-valued scalarly integrable  function is Pettis integrable. An analogue result, for weakly integrable fuzzy mappings has been proven in  \cite{Xue1} under the hypothesis that the target Banach space  is separable. In Theorem \ref{t4} and Theorem \ref{tcore2}  we present non-separable versions of  the above result.

\section{Basic facts.}

 Let $(\Omega,  {\mathcal L}, \mu)$ be a measure  space, where
 ${\mathcal L}$ is the family of all $\mu$-measurable sets.
Let $X$ be a Banach space endowed with the
norm $\|\cdot\|$ and let $X^*$ be its dual, $B(X)$ and $B(X^*)$ will denote respectively the unit ball in $X$ and in $X^*$. The symbol $c(X)$ stands for the family of all nonempty closed convex subsets of $X$, while $cwk(X)$ denotes the subfamily consisting of  weakly compact sets. In  $cwk(X)$ we consider the Minkowski addition ($A+B:=\{x+y: \ x\in A, y \in B\}$) and the
   standard multiplication by scalars ($k A :=\{kx: \ x\in A\}$). If $A \in cwk(X)$ and $x \in X$, by the symbol $A - x $ we denote the set $\{y \in X  :  y = z-x, \ z \in A\}$. Moreover we endow $cwk(X)$ with the {\it Hausdorff distance}

$$d_H(A,B):=\max\left\{\sup_{x \in A} \inf_{y \in B}\|x-y\|, \sup_{y
\in B}
\inf_{x \in A}\|x-y\|\right\} .
$$

\noindent The space $cwk(X)$   with the Hausdorff
   distance is a complete metric space. For every $A \in cwk(X)$ the {\it support
  function }of  $A$ is denoted by $s( \cdot, A)$ and defined  by
$s(x^*, A) := \sup \{ \langle x^*,y \rangle : \ y
  \in A\}$, for each $x^* \in X^*$. Clearly the map $x^*
\longmapsto s(x^*, A)$ is sublinear on $X$ and $-s(-x^*, A) =
\inf \{ \langle x^*,y \rangle : \ y \in A\}$, for each $x^* \in
X^*$.

\noindent Below, by a series $\sum_{n=1}^{\infty} A_n$, where  $A_n \in 2^X \setminus \{\emptyset\}$ for each $n$, we understand the set
$$\left\{\sum_{n=1}^{\infty} x_n: \ x_n\in A_n \ \ {\rm and}  \ \ \sum_{n=1}^{\infty} x_n \ \ {\rm unconditionally \ convergent}\right\}.$$

\noindent Let $u\colon X \to [0, 1]$. We set $[u]^r = \{x \in X: \  u(x) \geq r\}$, for $r\in (0,1]$;  $u$ is called a  {\it generalized fuzzy number}, as in   \cite{Xue1} and \cite{Xue},  ({\it fuzzy number} as in \cite{WW}) on  $X$ if, for each $r\in (0,1]$, $[u]^r \in {cwk(X)}$. Let ${\mathcal F}(X)$ denote the {\it generalized fuzzy number space}.

\noindent We define $\theta\colon X \to [0, 1]$ as follows:

$$ \theta(x) =\left \{ \begin{array}{lllll}
 1   & if  &x=0\\
 0   & if  &x \not=0 .
\end{array}
\right . $$

\noindent Then $ \theta \in {\mathcal F}(X)$ and $ \theta$ is called the {\it null element} of ${\mathcal F}(X)$.

\noindent In the sequel we will use the following representation theorem (see
 \cite{WW}).

\begin{thm}\label{t1} If $u \in {\mathcal F}(X)$, then
\begin{itemize}
\item[$(i)$] $[u]^r \in cwk(X)$, for all $r \in (0, 1]$;
\item[$(ii)$] $[u]^{r_2} \subseteq [u]^{r_1}$, for $0 \leq r_1 \leq r_2 \leq
1$;
\item[$(iii)$] if $(r_k)$ is a nondecreasing sequence converging to $r > 0$,
then

$$[u]^r= \bigcap_{ k \geq 1} [u]^{r_k}.$$
\end{itemize}
Conversely, if $\{A_r\colon \ r\in (0,1]\}$ is a family of subsets
of $X$ satisfying $(i)$--$(iii)$, then there exists a unique $u \in {\mathcal F}(X)$ such
that $[u]^r = A_r$ for $r\in (0,1]$.
\end{thm}

\noindent For $u, \ v\in {\mathcal F}(X)$ and $k \in \R$, define $u+ v$ and $\lambda u$ as follows (see \cite[pag. 536]{Park1})
\bee(u+v)(x)= \sup_{x= y+z} \min (u(y), v(z)),\eee
\bee (\lambda u)(x) = \left \{ \begin{array}{lllllll}
u \left(\frac{1}{\lambda} x \right)   & if  & \lambda  \not= 0  & \\
\\
\wt{0} & if &  \lambda  = 0, & {\rm where} \ \   \wt{0}= \chi_{\{0\}}
\end{array}
\right . \eee
\noindent For $u, \ v\in {\mathcal F}(X)$ and $k \in \R$, 
the addition and the scalar
multiplication are defined respectively by 
$$[u+v]^r:=[u]^r+[v]^r \ \ \mbox{and} \ \ [ku]^r:=k[u]^r ,$$

\noindent for each $r \in(0,1]$.

\begin{deff}  {\rm (see \cite[Definition 3.1]{Xue}).   A mapping $\wt{M}: {\mathcal L} \to {\mathcal F}(X)$ is called a {\it generalized fuzzy number measure}, if it satisfies the following requirements:

\begin{enumerate}
\item[$(1)$] $\wt{M(}\phi)= \theta$;
\item[$(2)$] $\wt{M}(\cdot)$ is countably additive, i.e. for any sequence $(A_n)_{n=1}^{\infty}$ of pairwise disjoint elements of $ {\mathcal L}$ we have that $\wt{M}(\cup_{n=1}^{\infty}A_n) = \sum_{n=1}^{\infty}\wt{M}(A_n)$, where $\sum_{n=1}^{\infty}\wt{M}(A_n)$ is defined by

$$ \sum_{n=1}^{\infty}\wt{M}(A_n)(x) = \sup  \ \{\inf_n \wt{M}(A_n)(x_n) :  \ \  x= \sum_{n=1}^{\infty} x_n  \ \  {\rm (unconditionally \ \ convergent)}\}.$$
\end{enumerate}}
\end{deff}

\noindent We recall that $M: {\mathcal L} \to 2^X \setminus \{\emptyset\}$ is a {\it weak  multimeasure}   if  $s(x^*, M)$ is a real valued measure for all $x^* \in X^*$,   $M: {\mathcal L} \to 2^X  \setminus \{\emptyset\}$    is a {\it set valued measure} if it satisfies the following two requirements
\begin{itemize}
\item[$(i)$] $M(\phi)= \{0\}$;
\item[$(ii)$] $M(\cdot)$ is countably additive, i.e. for any sequence $(A_n)_{n=1}^{\infty}$ of pairwise disjoint elements of $ {\mathcal L}$ we have that $M(\cup_{n=1}^{\infty}A_n) = \sum_{n=1}^{\infty}M(A_n)$.
\end{itemize}

\noindent Each set-valued measure is a  weak multimeasure.  If $M$ is a weak multimeasure taking values in $cwk(X)$ also the reverse implication holds (see \cite[Proposition 7.4]{hp}).

In the sequel we will use the following representation theorem (see
\cite[Lemma 3.3]{WuWu}).

\begin{thm}\label{t2} Let $\wt{M}: {\mathcal L} \to {\mathcal F}(X)$ be a  generalized fuzzy number measure, $M^{r}(A)= \{x \in X : \ \wt{M}(A)(x) \geq r\}$, $A \in {\mathcal L}$, $r \in (0,1]$. Then

\begin{itemize}

\item[$(i)$] $M^{r}: {\mathcal L} \to cwk(X)$  is a set valued measure, $r \in (0,1]$;
\item[$(ii)$] $M^{r_2}(A) \subseteq M^{r_1}(A)$, for $0 \leq r_1 \leq r_2 \leq
1$ and for  each  $A \in {\mathcal L}$;
\item[$(iii)$] if $(r_k)$ is a nondecreasing sequence converging to $r > 0$ and $A \in {\mathcal L}$,
then
$$M^r(A)= \bigcap_{ k \geq 1} M^{r_k}(A).$$
\end{itemize}
Conversely, if there is a family $\{M^r : \  r \in (0,1]\}$ satisfying the requirements (i)-(iii), we define a mapping $\wt{M}:  {\mathcal L} \to  {\mathcal F}(X)$ as

\be\label{AA} [\wt{M}(A)](x) = \left \{\begin{array}{lll}
 \sup\{r: x \in M^r(A), \ r \in (0,1]\}   & if  &\{r: x \in M^r(A), \ r \in (0,1]\}  \not= \phi,\\
 0   & if & \{r: x \in M^r(A), \ r \in (0,1]\} = \phi
\end{array}
\right . \ee
Then, $\wt{M}$ is a generalized fuzzy number measure and $[\wt{M}(A)]^r= M^r(A)$ for each  $A \in {\mathcal L}$ and  $r \in (0,1]$.
\end{thm}

We recall now some definitions concerning set-valued functions.
\noindent A set-valued function $\vG  \colon \Omega \to cwk(X)$  is said to be {\it scalarly measurable} if for
every $ x^* \in X^*$, the map
  $s(x^*,\vG)$ is measurable.  A set-valued function $\vG  \colon \Omega \to cwk(X)$ is
said to be {\it  scalarly}  {\it integrable on} $\Omega$ if for each $x^* \in
X^*$ the real function $s(x^*, \vG)$ is integrable  on $\Omega$.\\

\begin{deff}{\rm  (see \cite[Definition 3.1]{eh}) A set-valued function $\vG \colon \Omega \to
cwk(X)$ is said to be  {\it Pettis integrable}  on $\Omega$  if $\vG$ is
    scalarly  integrable on $\Omega$ and for each $A \in \mathcal{L}$
   there exists a set $W_{\Gamma}(A) \in cwk(X)$ such that for each
   $x^* \in X^*$, we have

  $$
   s(x^*,W_{\Gamma}(A))=\int_As(x^*,\vG) \,d \mu .$$
   Then we set $(P)\int_A\vG \,d \mu=: W_{\Gamma}(A) $, for each $A \in \mathcal{L}$.}

\end{deff}
\noindent Observe that $W_{\vG}$ is a $\mu$-continuous weak multimeasure.

\noindent A function  $f \colon \Omega \to X$ is called a {\it selection} of a set-valued
function $\vG  \colon\Omega \to cwk(X)$ if,
for every $\omega \in \Omega$, one has $f(\omega) \in\vG(\omega)$.
If $\vG \colon \Omega \to cwk(X)$ is set-valued function, by the symbol
${\mathcal S}(\vG)$ we denote the family of all Pettis integrable selections of $\vG$.
From now on if $\wt{\vG}\colon \Omega \to
{\mathcal F}(X)$  is a  fuzzy mapping, we set
$$\wt{\vG}_{r}(\omega):=[\wt{\vG}(\omega)]^r.$$
Given a vector valued function $f\colon \Omega \to X$ we call {\it fuzzy-number valued function associated to} $f$ the fuzzy  mapping $\wt{f}\colon \Omega \to {\mathcal F}(X)$ defined by $\wt{f}(\omega):=\chi_{f(\omega)}$. Using the same terminology as for multifunctions, if a function   $f \colon \Omega\to X$ is such that $f(\omega) \in
\wt{\vG}_{r}(\omega)$, for each $r \in (0,1]$, we say that the associated fuzzy mapping $\wt{f}(\omega)$ is a selection of $\wt{\vG}$.
We recall that a  fuzzy mapping $\wt{\vG}$ is said to be {\it scalarly} {\it integrable} on $\Omega$ if for all $r\in
(0,1]$ the set-valued function $\wt{\vG}_{r}\colon  \Omega \to cwk(X)$ is scalarly  integrable.

\section{Fuzzy Pettis integral}

We start with the definition of fuzzy Pettis integral.

\begin{deff} {\rm A fuzzy mapping $\wt{\vG}\colon\Omega \to{\mathcal F}(X)$
is said to be {\it  weakly Pettis   integrable }on $\Omega$ if for every
$r\in (0,1]$ the set-valued function $\wt{\vG}_{r}$ is
Pettis integrable on $\Omega$.}
\end{deff}

\begin{deff} {\rm A fuzzy mapping $\wt{\vG}\colon \Omega \to{\mathcal F}(X)$
is said to be  {\it  Pettis   integrable} on $\Omega$ if it is weakly Pettis integrable  and,  for each $A \in {\mathcal L}$, there exists a generalized fuzzy number
$\wt{M_{\vG}}(A)  \in{\mathcal F}(X)$ such that for any  $r\in(0,1]$ and for any $x^*\in X^*$ we have

\begin{equation}\label{10}s(x^*,[\wt{M_{\vG}}(A)]^r)=\int_A s(x^*,\wt{\vG}_{r})\,d \mu,
\end{equation}

\noindent (see \cite[Definition 3.2]{Park1}).

\noindent We call the fuzzy mapping $\wt{M_{\vG}}\colon {\mathcal L}\to{\mathcal F}(X)$  the {\it fuzzy Pettis integral of} $\wt{\vG}$ and we set
$(P)\int_A\wt{\vG}\,d \mu :=\wt{M_{\vG}}(A)$.}
\end{deff}

\noindent It follows at once from the definition that a vector valued function $f:\Omega \to X$ is Pettis integrable if and only if the associated fuzzy mapping  $\wt{f}\colon \Omega \to {\mathcal F}(X)$ is so.
We recall the following properties:
\begin{prop}\label{prp0} {\rm (see \cite[Theorem 3.3]{Park1}) \ \ Let  $\wt{F}\colon \Omega \to{\mathcal F}(X)$ and   $\wt{G}\colon \Omega\to{\mathcal F}(X)$ be  Pettis integrable and let $\lambda \geq 0$. Then

\begin{enumerate}
\item[$(1)$] $\wt{F} +\wt{G}$ is Pettis integrable and  for each $A \in {\mathcal L}$

$$(P)\int_A (\wt{F} +\wt{G})\,d \mu= (P)\int_A \wt{F}\,d \mu + (P)\int_A \wt{G}\,d \mu
;$$
\item[$(2)$] $\lambda\wt{F}$ is Pettis integrable and  for each $A \in {\mathcal L}$

$$(P)\int_A \lambda\wt{F} \,d \mu= \lambda  \ (P)\int_A \wt{F}\,d \mu .$$
\end{enumerate}
}\end{prop}
\noindent The following decomposition result  will be useful in our investigation.

\begin{thm}\label{tdec} Let  $\wt{\vG}\colon \Omega \to{\mathcal F}(X)$ be a fuzzy
 scalarly integrable mapping. Then, the following conditions are equivalent:
\begin{enumerate}
\item[$(i)$]  $\wt{\vG}$ is fuzzy Pettis integrable on $\Omega$;
\item[$(ii)$] ${\mathcal S}(\wt\vG_1) \not= \emptyset$ and
for every  $f\in{\mathcal
S}(\wt\vG_1)$,   there exists a  fuzzy  Pettis integrable mapping  $\wt{G}\colon \Omega \to{\mathcal F}(X)$ such that  $\wt{\vG}(\omega)=\wt{G}(\omega)+\wt{f}(\omega)$ (where $\wt{f}$ is the fuzzy mapping associated to $f$)  and,  for each $r \in (0,1]$ and for each $\omega \in\Omega$, $0 \in \wt{G}_r(\omega)$.
\end{enumerate} \end{thm}

\noindent {\Pf }  If $\wt{\vG}$ is a  Pettis integrable fuzzy mapping on $\Omega$, for each $r \in (0,1]$ the multifunction $ \wt{\vG}_{r}$ is Pettis integrable. Therefore ${\mathcal S}(\wt{\vG_1}) \not= \emptyset$ (see \cite{cascales kadets rodriguez1} or \cite{mu}).

\noindent Let us fix $f\in{\mathcal S}(\wt\vG_1)$ and  let $\wt{f}$ be  the fuzzy-number valued function associated to $f$.  Now  define $\wt{G}\colon \Omega \to{\mathcal F}(X)$  setting
$\wt{G}_r(\omega):=\wt{\vG}_r(\omega)-f(\omega)$ for each  $r \in (0,1]$. At first we observe  that $\wt{G}$ is a fuzzy mapping. In fact,  since the sets $\wt{G}_r(\omega)$ are translations of the
 sets $\wt{\vG}_r(\omega)$ with a same vector,    conditions $(i)-(iii)$ of Theorem \ref{t1} are trivially satisfied.

\noindent  We are going  to prove that $\wt{G}$ is   Pettis integrable on $\Omega$. Since for each  $r \in (0,1]$,  $\wt{\vG}_r$ is Pettis  integrable and $f$ is Pettis integrable,  then also $\wt{G}_r$ is Pettis integrable for each  $r \in (0,1]$.

\noindent  Moreover denote by  $F$ the Pettis  primitive of $f$ and by $\wt{M_{\vG}}$ the fuzzy Pettis integral of $\wt{\vG}$.  Now  for each  $A \in {\mathcal L}$ define $\wt{N_{G}}(A)$  setting $$[\wt{N_{G}}(A)]^r:=[\wt{M_{\vG}}(A)]^r-F(A)$$  for  each  $r \in (0,1]$.

\noindent  Also in this case we have that  $\wt{N_{G}}(A)$ is a fuzzy number. Moreover for each $x^* \in X^*$ and for  each  $r \in (0,1]$ we have
\bee \begin{split} s(x^*,[\wt{N_{G}}(A)]^r) &= s(x^*, [\wt{M_{\vG}}(A)]^r) -x^*(F(A)) =\int_A s(x^*,\wt{\vG}_{r})d \mu   -\int_A x^*fd \mu \\
&=\int_A s(x^*, \wt{\vG}_r-f)d \mu= \int_A s(x^*, \wt{G}_r)d \mu.
\end{split}\eee
\noindent Therefore $\wt{G}$ is Pettis is integrable. Since     for each $r \in (0,1]$ and for each $\omega \in \Omega$, $\wt{\vG}_{1}(\omega)\subset \wt{\vG}_{r}(\omega)$,  then $f\in{\mathcal S}(\wt\vG_r)$. Therefore  $0 \in \wt{G}_r(\omega)$.

\noindent $(ii) \Longrightarrow (i)$ follows from $(1)$ of Proposition (\ref{prp0}).

{\f }

\vspace{2ex}

\noindent Observe that $0 \in \wt{G}_r(\omega)$ implies that the support function $s(x^*,\wt{G}_r)$ is nonnegative.

\noindent We are going to prove that the fuzzy Pettis integral is a generalized fuzzy number measure.

\begin{thm}\label{p1} Let   $\wt{\vG}\colon \Omega \to {\mathcal F}(X)$ be
    a Pettis integrable
fuzzy mapping on $\Omega$. Then the fuzzy Pettis integral  $\wt{M_{\vG}}$ of $\wt{\vG}$  is a generalized fuzzy number measure.
    \end{thm}

\noindent {\bf Proof.}  If  $\wt{\vG}$  is a Pettis integrable fuzzy mapping, its Pettis integral is a fuzzy mapping $\wt{M_{\vG}}\colon {\mathcal L}\to{\mathcal F}(X)$. Moreover
 for each $r \in (0,1]$, the multifunction $\wt{\vG}_{r}: \Omega \to cwk(X)$ is Pettis integrable and, for each $A \in {\mathcal L}$,

 $$[\wt{M_{\vG}}(A)]^r=(P)\!\int_A\wt{\vG}_{r}\,d \mu.$$

\noindent  Therefore $[\wt{M_{\vG}}]^r$  is a weak multimeasure.  Since  $[\wt{M_{\vG}}]^r$  takes values in $cwk(X)$,  it is a multimeasure, and property (i) of Theorem \ref{t2} is satisfied.

\noindent Moreover properties (ii) and (iii) of Theorem \ref{t2} are satisfied since, for each $A \in {\mathcal L}$,  $\wt{M_{\vG}}(A)$ is a fuzzy number. Therefore the family  $\{[\wt{M_{\vG}}(A)]^r : \  r \in (0,1]\}$ satisfies condition (i)--(iii) of Theorem \ref{t2}. Then
the mapping $\wt{M}:  {\mathcal L} \to  {\mathcal F}(X)$ defined as in (\ref{AA}) is a fuzzy number measure with $\wt{M}^r(A)= [\wt{M_{\vG}}(A)]^r$ for each  $A \in {\mathcal L}$ and  $r \in (0,1]$. Moreover $\wt{M}(A)= \wt{M_{\vG}}(A)$ and the assertion follows.
 \hfill$\Box$

\vspace{2ex}

\noindent Let   $\wt{\vG}\colon \Omega \to {\mathcal F}(X)$ be a Pettis integrable fuzzy mapping. Then  for each $r \in (0,1]$,
$\wt{\vG}_{r}: \Omega \to cwk(X)$ is a Pettis integrable multifunction. Therefore (see \cite[Theorem 1.4]{mu}), for each $r \in (0,1]$, the operator
 $T_r: X^* \to L^1(\mu)$, defined as $T_r(x^*) := s(x^*, \wt{\vG}_{r})$ is $\tau(X^*, X)$- weakly continuous, where $\tau(X^*, X)$ denotes the topology of uniform convergence on weakly compact subsets of $X$.

 \noindent In case of fuzzy mappings, we have the following characterization.

\begin{thm}\label{t3}   Let   $\wt{\vG} \colon\Omega \to {\mathcal F}(X)$ be
    a scalarly integrable  fuzzy mapping. Then
     $\wt{\vG}$  is Pettis integrable if and only if for  each   $r \in (0,1]$ the operator
 $T_r$ is $\tau (X^*, X)$-weakly continuous. If   $0 \in \wt{\vG}_r(\omega)$  for each $r \in (0,1]$  and for each $\omega \in \Omega $, then $T_r$ is $\tau(X^*, X)$-norm continuous.

 \end{thm}
\noindent {\bf Proof.}
Let  $\wt{\vG}\colon \Omega\to {\mathcal F}(X)$ be a Pettis integrable fuzzy mapping. Then  for each $r \in (0,1]$,
 $\wt{\vG}_{r}:\Omega \to cwk(X)$ is a Pettis integrable multifunction therefore by \cite[Theorem 1.4]{mu}  the operator
 $T_r$ is $\tau(X^*, X)$-weakly continuous.

\noindent Conversely assume
 that for each $r \in (0,1]$ the operator
 $T_r$ is $\tau(X^*, X)$-weakly continuous. Again by  \cite[Theorem 1.4]{mu}  it follows that $\wt{\vG}_{r}:\Omega \to cwk(X)$ is a Pettis integrable multifunction.
 For each $A \in {\mathcal L}$ and for each $r \in (0,1]$ define $[\wt{M_{\vG}}(A)]^r \in cwk(X)$   such that
\bee
s(x^*,[\wt{M_{\vG}}(A)]^r)=\int_As(x^*,\wt{\vG}_{r})\,d \mu.
\eee

 \noindent 
We have  to prove that for any   $A \in {\mathcal L}$
the family  $\{ [\wt{M_{\vG}}(A)]^r : r \in (0,1] \}$,
satisfies  properties (i)--(iii) of Theorem \ref{t1}.
Clearly condition (i) is satisfied.
In order to prove  (ii), let $0 \leq r_1
\leq r_2 \leq 1$. Since  $\wt{\vG}$ is a fuzzy mapping,  by Theorem \ref{t1} we have
$\wt{\vG}_{r_2}(\omega)\subseteq \wt{\vG}_{r_1}(\omega)$, for each $\omega \in \Omega$. Therefore
$$s(x^*,\wt{\vG}_{r_2}) \leq s(x^*, \wt{\vG}_{r_1}),$$
for each $x^* \in X^*$. Therefore for each $A \in {\mathcal L}$  and for each $x^* \in X^*$

\be\label{ec}s(x^*, [\wt{M_{\vG}}(A)]^{r_2})= \int_A s(x^*,\wt{\vG}_{r_2})\,d \mu \leq \int_A s(x^*, \wt{\vG}_{r_1})\,d \mu= s(x^*,[\wt{M_{\vG}}(A)]^{r_1}).\ee

\noindent Then, as a consequence of  Hahn-Banach separation
theorem for convex sets, by (\ref{ec}) we also infer the inclusion      $[\wt{M_{\vG}}(A)]^{r_2}
\subseteq [\wt{M_{\vG}}(A)]^{r_1}$  for each $A \in {\mathcal L}$  and property (ii) is satisfied.
If  $(r_k)$ is
a nondecreasing sequence converging to $r > 0$, then for each $\omega \in
\Omega$ we have

$$\wt{\vG}_r(\omega)= \bigcap_{ k \geq 1} \wt{\vG}_{r_k}(\omega)\,.$$
\noindent Consequently (see \cite[Proposition 1]{so})
$$s(x^*, \wt{\vG}_r)= \lim_k s(x^*, \wt{\vG}_{r_k})\,,$$
for each  $\omega \in \Omega$ and $x^* \in X^*$.

\noindent By hypothesis, for each  $x^* \in X$, the sequence of real
valued functions $\left(s(x^*, \wt{\vG}_{r_k})\right)$ is integrable.
As in Theorem \ref{tdec}, we can observe that, for each $r \in (0,1]$ and for each  $\omega \in \Omega$,  $0 \in \wt{G}_r(\omega) = \wt{\vG}_{r}(\omega) -  f(\omega)$. Then, without loss of generality, we can assume that $s(x^*, \wt{\vG}_{r_k}) \geq 0$. Since  $s(x^*, \wt{\vG}_{r_1}) \geq s(x^*, \wt{\vG}_{r_k}\geq 0$, by Lebesgue dominated convergence theorem we get

\begin{eqnarray*}
 \  s(x^*,  [\wt{M_{\vG}}(A)]^r)&=& \int_A s(x^*,
\wt{\vG}_r)\,d \mu \
\\
&=&\lim_k  \int_A s(x^*, \wt{\vG}_{r_k})\,d \mu\\
&=&\lim_k  s(x^*,  [\wt{M_{\vG}}(A)]^{r_k})= s(x^*, \bigcap_{ k \geq 1} [\wt{M_{\vG}}(A)]^{r_k})
.\end{eqnarray*}

\noindent Since  above equalities  hold for each  $x^* \in X^*$,  we obtain
$ [\wt{M_{\vG}}(A)]^r=\bigcap_{ k \geq 1} [\wt{M_{\vG}}(A)]^{r_k}$ and property (iii) is satisfied.

\noindent If $0 \in \wt{\vG}_r(\omega)$  for each $r \in (0,1]$ then again by \cite[Theorem 1.4]{mu} it follows that $T_r$ is $\tau(X^*, X)$-norm continuous.

\hfill$\Box$

 \vspace{2ex}

\noindent  We recall that $X^*$ is said to be {\it weak$^*$- angelic} if for each bounded set $B \subset X^*$ the weak$^*$-closure of $B$ is equal to the set of weak$^*$-limits of sequences from $B$.

\begin{deff}  {\rm We say that a space $Y \subset X$ {\it determines a  fuzzy mapping} $\wt{\vG}\colon \Omega\to {\mathcal F}(X)$, if for each $r \in (0,1]$, $s(x^*, \wt{\vG}_{r}) = 0$ $ \mu-a.e.$ for each $x^* \in Y^{\bot}$, where $Y^{\bot}$ denotes the annihilator of $Y$ in $X^*$ and the exceptional sets depend on $x^*$ and on $r$}.
\end{deff}

 \begin{thm}\label{t4} Let    $\wt{\vG}\colon \Omega\to {\mathcal F}(X)$ be a scalarly integrable fuzzy mapping determined by a space $Y \subseteq X$ such that $Y^*$ is weak$^*$-angelic and such that for each $r \in (0,1]$ the operator  $T_r: X^* \to L^1(\mu)$ is weakly compact. Then $\wt{\vG}$  is a Pettis integrable  fuzzy mapping.
    \end{thm}

\noindent {\bf Proof.}
Let $\wt{\vG}\colon \Omega\to {\mathcal F}(X)$ be  a scalarly integrable fuzzy mapping satisfying the hypothesis of the claim. Then  for each $r \in (0,1]$,
$\wt{\vG}_{r}:\Omega \to cwk(X)$ is a scalarly integrable multifunction and satisfies hypotheses (WC) and (D$^*$) of  \cite[Theorem 2.4]{mu}.  Therefore  it follows that $\wt{\vG}_{r}$ is a Pettis integrable multifunction. For each $A \in {\mathcal L}$ and for each $r \in (0,1]$ let $[\wt{M_{\vG}}(A)]^r \in cwk(X)$  be such that equality (\ref{10}) holds. 

\noindent We have to prove that for any   $A \in {\mathcal L}$
the family  $\{[\wt{M_{\vG}}(A)]^r : \  r \in (0,1]\}$
satisfies  properties (i)--(iii) of Theorem \ref{t1}.  To this end it is enough to proceed as in the proof of the $ `` only \  if ''$
part of  Theorem \ref{t3}. Thus $\wt{M_{\Gamma}}$ defines a fuzzy mapping and we conclude that $\wt{\vG}$ is  fuzzy Pettis integrable.       \hfill$\Box$

\begin{deff}  {\rm A subspace $Y \subset X$ is said to be {\it weakly compact generated} (WCG) if there is a weakly compact set $K$ in $Y$ such  that $Y=\overline{ \rm span}(K)$.}
\end{deff}

\noindent Examples of WCG spaces are separable spaces, reflexive spaces, $L^1(\mu)$ if $\mu$ is a $\sigma$-finite measure, and $c_0(\Gamma)$ if $\Gamma$ is a nonempty set (see \cite[Chapters 13]{Fabian1}).

 \begin{prop}\label{p2} Let   $\wt{\vG}\colon \Omega\to {\mathcal F}(X)$ be a Pettis integrable fuzzy mapping. Then $\wt{\vG}$ is determined by a WCG subspace of $X$  and  for each $r \in (0,1]$ the operator  $T_r: X^* \to L^1(\mu)$ is weakly compact.
 \end{prop}

\noindent  {\bf Proof.}  If  $\wt{\vG} \colon \Omega\to {\mathcal F}(X)$ is Pettis integrable, then  for each $r \in (0,1]$,
$\wt{\vG}_{r}:\Omega \to cwk(X)$ is a Pettis integrable multifunction. Therefore by \cite[Proposition 2.2]{mu} the operator  $T_r: X^* \to L^1(\mu)$ is weakly compact. Moreover as in  \cite[Proposition 2.2]{mu},  denote by $Y_r$ the weakly compact generated subspace of $X$ which determines the multifunction $\wt{\vG}_{r}$.
If $0 < r_1 \leq r_2 \leq 1$, then $\wt{\vG}_{r_2}(\omega) \subseteq \wt{\vG}_{r_1}(\omega) $, and also $Y_{r_2} \subseteq Y_{r_1} $.   Let now $\{s_k\}$ be a sequence in $(0,1]$ such that $s_k \searrow 0$ and for each $k $ let $W_{s_k} \subset B(X)$ be a weakly compact convex set generating $Y_{s_k}$. The set $\bigcup_{k} \frac{1}{2^k}W_{s_k}$ is a weakly compact set generating a space $Y$. So if $x^* \in Y^{\bot}$,  for each $s_k$,  $x^* \in Y_{s_k}^{\bot}$ and then also $x^* \in Y_{r}^{\bot}$  for each $r \in (0,1]$. It follows from   \cite[Proposition 2.2]{mu}  that $s(x^*, \wt{\vG}_{r}) = 0$ $ \mu-a.e.$  This ends the proof.
 \hfill$\Box$

\vspace{2ex}

\noindent We recall that the class of Banach spaces having ${\rm weak}^*$-angelic dual is very large and contains all weakly compactly generated spaces (\cite[Chapters 11-12]{Fabian}). Then the following theorem is a consequence the previous Theorem \ref{t4} and Proposition \ref{p2}.

\begin{thm}\label{t5} A scalarly integrable fuzzy mapping $\wt{\vG}\colon \Omega\to {\mathcal F}(X)$ is Pettis integrable if  and only if  $\wt{\vG}$ is determined by a WCG space $Y \subseteq X$ and  for each $r \in (0,1]$ the operator  $T_r: X^* \to L^1(\mu)$ is weakly compact.
\end{thm}

\noindent As a consequence of Theorem \ref{t5} we obtain:

\begin{cor}\label{c1} Let  $\wt{\vG}\colon \Omega\to {\mathcal F}(X)$ be a scalarly measurable  fuzzy mapping. Assume that there exists a fuzzy mapping   $\wt{G}\colon \Omega\to {\mathcal F}(X)$ which is  Pettis integrable and satisfies, for each $x^* \in X^*$ and $r \in (0,1]$,

\be\label{c1.1} | s(x^*,  \wt{\vG}_{r})| \leq  | s(x^*,  \wt{G}_{r})| .\ee
\noindent Then  $\wt{\vG}$ is Pettis integrable.
\end{cor}
\noindent {\Pf} Indeed inequality (\ref{c1.1}) implies that  for each $r \in (0,1]$, $\wt{\vG}$ is determined by a WCG space $Y \subseteq X$  and the weak compactness of the operator  $T_r: X^* \to L^1(\mu)$, for each $r \in (0,1]$.
{\f}

\vspace{2ex}

\noindent It is known that if $X$ is a separable Banach space and $c_0 \not\subseteq X$, then each scalarly integrable function is Pettis integrable. An analogue result, for weakly integrable bounded fuzzy mapping has been proved in  \cite[Theorem 4.5]{Xue1} in case of a separable space. In the following theorem  we prove  Pettis integrability of a scalarly integrable fuzzy mapping  in case of a general Banach space without any boundness condition, extending  \cite[Theorem 2.13]{mu}   proved for multifunctions.

\begin{thm}\label{t6} Let  $X$ be a Banach space not containing any isomorphic copy of $c_0 $.  If  $\wt{\vG}\colon \Omega\to {\mathcal F}(X)$ is a scalarly measurable  fuzzy mapping which  is determined by a WCG space, then $\wt{\vG}$ is Pettis integrable.
\end{thm}
\noindent {\Pf} If  $\wt{\vG}\colon \Omega\to {\mathcal F}(X)$ is a scalarly integrable fuzzy mapping, then  for each $r \in (0,1]$,
$\wt{\vG}_{r}:\Omega \to cwk(X)$ is a scalarly integrable multifunction and is determined by a WCG space. Therefore by  \cite[Theorem 2.13]{mu} it follows that $\wt{\vG}_{r}:\Omega \to cwk(X)$ is a Pettis integrable multifunction.

\noindent Let  $[\wt{M_{\vG}}]^r $ be the Pettis integral of  $\wt{\vG}_{r}$.  It is enough to prove that for each $A \in {\mathcal L}$ the family  $\{ [\wt{M_{\vG}}(A)]^r : r \in (0,1] \}$,
satisfies  properties (i)--(iii) of Theorem \ref{t1}.  So, proceeding as in  the proof of the $ `` only \  if ''$
part of the proof Theorem \ref{t3}, we get the result.

{\f}

\vspace{2ex}

 \section{Core characterization}

The result of Talagrand which characterizes Pettis integrability of a vector valued function $f$ is extended to multifunction in \cite{mu}. We recall that (see \cite{mu} Definition 4.3) if $\vG  \colon \Omega \to cwk(X)$ is a multifunction, for  each $E \in \Sigma$,  the {\it core} of $\vG$ {\it on $E$} is defined by

 $${\rm cor}_{\vG }(E)= \bigcap_{\mu(N)=0} \overline{\rm conv} \vG(E \setminus N)=  \bigcap_{\mu(N)=0} \overline{\rm conv}\left( \bigcup_{\omega \in E \setminus N} \vG(\omega) \right).$$

 \noindent We say that a fuzzy mapping  $\wt{G} \colon \Omega\to {\mathcal F}(X)$ is {\it dominated by} a fuzzy mapping  $\wt{\vG}\colon \Omega\to {\mathcal F}(X)$ if, for  each $r \in (0,1]$ and for each $\omega \in \Omega$,  $\wt{G}_{r}(\omega) \subseteq \wt{\vG}_{r}(\omega)$.

\noindent  The following theorem extends to scalarly integrable fuzzy mappings the core characterization of Pettis integrability.

\noindent The symbol  $\Sigma^+$ denotes the sets  $E \in \Sigma$ with $\mu(E)>0$.

 \begin{thm}\label{tcore} Let   $\wt{\vG}\colon \Omega\to {\mathcal F}(X)$ be a scalarly integrable fuzzy mapping. Then $\wt{\vG}$ is Pettis integrable if and only if the following properties hold:
 \begin{enumerate}
\item[$(1)$] for each $r \in (0,1]$ the operator  $T_r: X^* \to L^1(\mu)$ is weakly compact;
\item[$(2)$] if $\wt{G}\colon \Omega\to {\mathcal F}(X)$ is a scalarly measurable fuzzy mapping that is dominated by $\wt{\vG}$, then for each $r \in (0,1]$, ${\rm cor}_{\wt{G}_{r}} (E) \not= \emptyset$, for every $E \in \Sigma^+$.
\end{enumerate}
 \end{thm}

\noindent {\bf Proof.}
Assume first that  $\wt{\vG}$ is Pettis integrable, then for each $r \in (0,1]$, the multifunction $\wt{\vG}_{r}:\Omega \to cwk(X)$ is Pettis integrable. Therefore the operator $T_r: X^* \to L^1(\mu)$ is weakly compact. Moreover if $\wt{G}\colon \Omega\to {\mathcal F}(X)$ is a scalarly measurable fuzzy mapping which is dominated by $\wt{\vG}$, we have that $ \wt{G}_r(\omega) \subseteq  \wt{\vG}_{r}(\omega)  $ for each $\omega$ and for each $r$.
Then applying \cite[Theorem 4.6]{mu}, we get that ${\rm cor}_{\wt{G}_{r}} (E) \not= \emptyset$, for every $E \in \Sigma^+$.

\noindent Conversely assume that conditions (1) and (2) are satisfied. This implies that for each  $r \in (0,1]$, the multifunction
$\wt{\vG}_{r}:\Omega \to cwk(X)$ is a scalarly integrable multifunction, such that  the operator $T_r: X^* \to L^1(\mu)$ is weakly compact. Moreover condition (2) implies that  
 ${\rm cor}_{\wt{G}_{r}} (E) \not= \emptyset$, for each $r$ and for every $E \in \Sigma^+$. 
Therefore, applying again  \cite[Theorem 4.6]{mu},   it follows that $\wt{\vG}_{r}:\Omega \to cwk(X)$ is a Pettis integrable multifunction.

\noindent If again  $[\wt{M_{\vG}}]^r $ is  the Pettis integral of  $\wt{\vG}_{r}$, then, by proceeding as in the proof of Theorem \ref{t3},  we infer that for each $A \in {\mathcal L}$ the family  $\{ [\wt{M_{\vG}}(A)]^r : r \in (0,1] \}$,
satisfies  properties (i)--(iii) of Theorem \ref{t1}. Therefore the assertion follows.

 \hfill$\Box$

 \noindent Observe that if $\wt{\vG}\colon \Omega\to {\mathcal F}(X)$ is a Pettis integrable fuzzy mapping, then each multifunction $\wt{\vG}_{r}:\Omega \to cwk(X)$ possesses Pettis integrable selections (see   \cite[Theorem 2.6]{cascales kadets rodriguez1}).  As a consequence of  previous theorem we obtain that

 \begin{thm}\label{tcore2} Let   $\wt{\vG}\colon \Omega\to {\mathcal F}(X)$ be a scalarly integrable fuzzy mapping. If for each $r \in (0,1]$ the operator  $T_r: X^* \to L^1(\mu)$ is weakly compact and every scalarly measurable selection of  $\wt{\vG}_{r}$ is Pettis integrable, then $\wt{\vG}$ is a Pettis integrable fuzzy mapping.

 \end{thm}

\noindent {\bf Proof.}  Let $\wt{G}\colon \Omega\to {\mathcal F}(X)$ be any scalarly measurable fuzzy mapping that is dominated by $\wt{\vG}$. Then, for each $r \in (0,1]$, the multifunction ${\wt{G}_{r}} $ has scalarly measurable selections (see \cite{ckr1}). It follows by the hypothesis that each scalarly measurable selection of ${\wt{G}_{r}} $ is Pettis integrable. Therefore ${\rm cor}_{\wt{G}_{r}} (E) \not= \emptyset$, for every $E \in \Sigma^+$. Consequently by Theorem \ref{tcore}, we get that
 $\wt{\vG}$ is a Pettis integrable fuzzy mapping.
 \hfill$\Box$

\vspace{4ex}

{\sc Acknowledgments.}  \ \ The authors are grateful to the anonymous reviewers for their valuable remarks.

\end{document}